# The Khipu-Based Numeration System

Sneha Pilgaonkar

**Abstract.** This paper presents a survey of the khipu-based enumeration system of the Incas. In this system, different knots are tied together and they are attached to a primary cord. The numerical basis of the khipu system is well understood but more recent theories suggest that the khipu also encoded narrative information and this may have been facilitated by binary mapping.

**Introduction**

The word "khipu" (or *quipu*) is a Quechua word for "knot". Quechua language was the language of Inca Empire during the period 1400-1532 CE whose capital was located in Cuzco City in high Andes of southern Peru; the language is still spoken in the region. The Empire encompassed the existing areas of entire Peru and some parts of Ecuador, Bolivia, Chile and Argentina. The khipu was used by the administrators of Inca Empire to record various types of data like census, tax payment, and so on. The people in the empire who made khipu were administrators and were known as *khipukamayuq* (knot-maker in Quechua).

Gary Urton of Harvard University has created a khipu database wherein data of all extant khipu are stored. A relational database is used for the project. The database consists of detailed information of all the types of cords in khipu- primary, pendant, top, subsidiary such as the material, method of spinning and plying, color etc. It also consists of details of the various knots, its directionality, its position on cord, its numerical value etc. [1]

The Inca were very proficient at mathematics and their enumeration system using knots is well understood [2],[3]. For general studies of the Inca Empire, see {4],[5]. Here's an early eye-witness account of the mathematical capabilities of the Inca:

> To see them use another kind of quipu, with maize kernels, is a perfect joy. In order to carry out a very difficult computation for which an able computer would require pen and paper, these Indians make use of their kernels. They place one here, three somewhere else and eight, I know not where. They move one kernel here and there and the fact is that they are able to complete their computation without making the smallest mistake. as a matter of fact, they are better at practical arithmetic than we are with pen and ink. Whether this is not ingenious and whether these people are wild animals let those judge who will! What I consider as certain is that in what they undertake to do they are superior to us. (Quoted by Kak [6])

The Spanish conquered the Inca Empire in the year 1532. The Spanish became aware of the many uses that the khipu was still being put to use. They caught an Indian worker who confessed that the khipu was being used to store information about the good and bad things about the ruler [7]. The Spanish destroyed many of the khipu but a few hundred still survive [1],[8]-[12]. To understand the khipu system, one may also see how numeration arose in other cultures [13-[22]. It should be mentioned that knotted strings were used as numerical devices in Greece, Persia, China, Ryukyu Islands, and even Germany [7]. But in none of these places knotted strings played their unique role in the Inca Empire. The connections between the khipu and astronomical computations are described in [23]-[27].

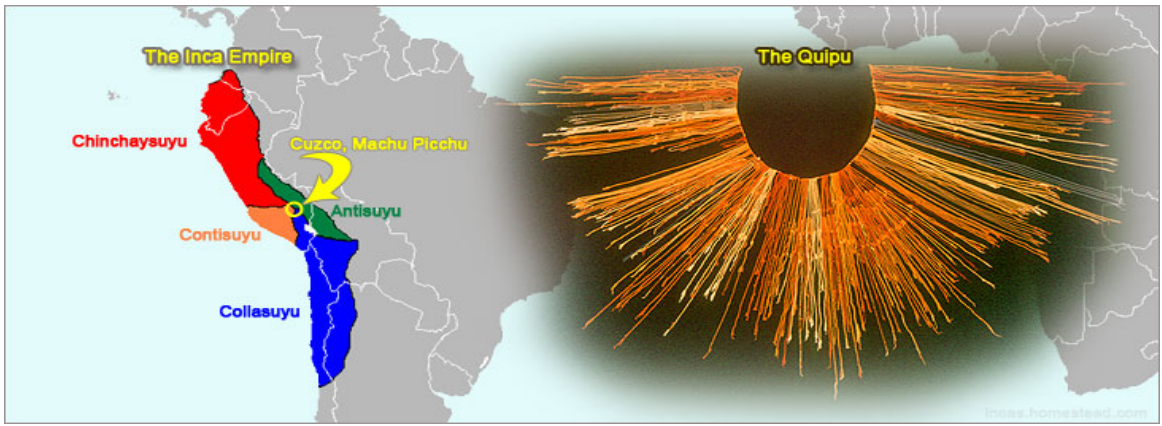

Fig.1 Provenance of the khipu [12]

This article presents a general overview of the khipu system. The problem of khipu coding is of interest not only for anthropological reasons, but also for abstract mathematical reasons related to information coding.

**Basic Construction of the Khipu**

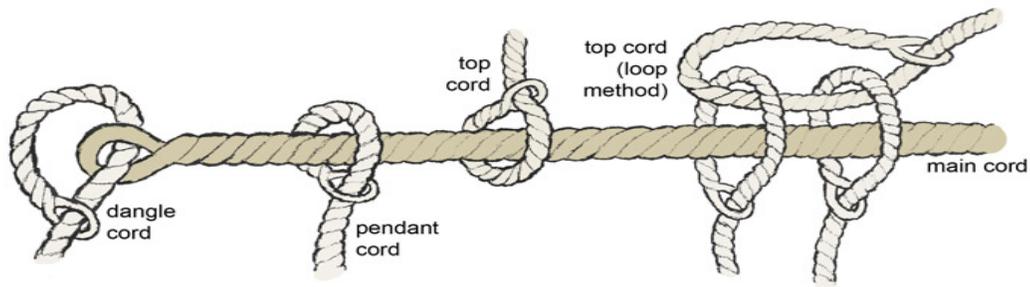

Fig.2 Basic structure of Khipu [10]

A khipu can be defined as several bunch of cryptic knotted strings made up of different material. Basically it is composed of several parts:

- Primary or Main chord- It the main chord of the khipu in either circular or horizontal position.
- Pendant chord- It is the chord which is attached to the primary chord and hands downwards. They are of different colors and different types of attachments.
- Subsidiary chords- These chords hang out from pendant chords. The khipu sometimes have 10-12 levels of subsidiary chords.
- Top chord- It is also a type of pendant chord attached to primary chord but in the upward direction opposite to the direction of pendant chord.
- Dangling chord – It is the chord which is at one end of main chord which is dangling hence the name dangling chord.

**Storing Information in Khipu**
Science historian Leland Locke was the first one to discover that around 100 khipu in American Museum in New York have encoded numerical information. Gary Urton proposes that khipu are unique since they are the only three dimensional 'written' documents all over the world. Further Urton also emphasizes that apart from numerical information, the khipu also stores narrative information, although not a single narrative khipu till date has been decoded. [9]

The khipu has encoded data in the following two formats:
1) Storing Numbers
2) Storing Narrative Information

**Storing Numbers in Khipu**
Numbers are stored in decimal format in a khipu. There are three types of knots which are used to encode decimal numbers in a khipu. These knots are located on pendants and its subsidiaries if any and the number value is added to get the value of the decimal number recorded on a khipu [8].

The long knot is used to record units place value of the number. Its value depends on how many turns of fibers are used. There is an exception for number '1' since long knot will have just one turn and it turns out to be same as simple knot. Figure Eight knot is used for storing number '1'. Single knot or simple knot is used for any places above units place i.e. hundreds, thousands etc. The number of turns on simple knot denotes its value in that decimal place. Units place is at the bottom of the pendant and as we go above towards the main chord, we encounter hundreds, thousands place etc. The highest value of the number i.e. the highest decimal place is the one nearest to the primary chord.

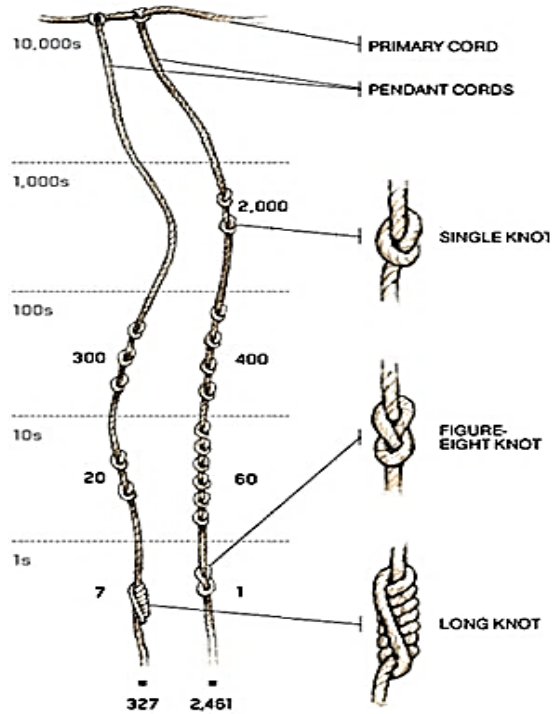

Fig.3 Reading Numbers in the Khipu [4]

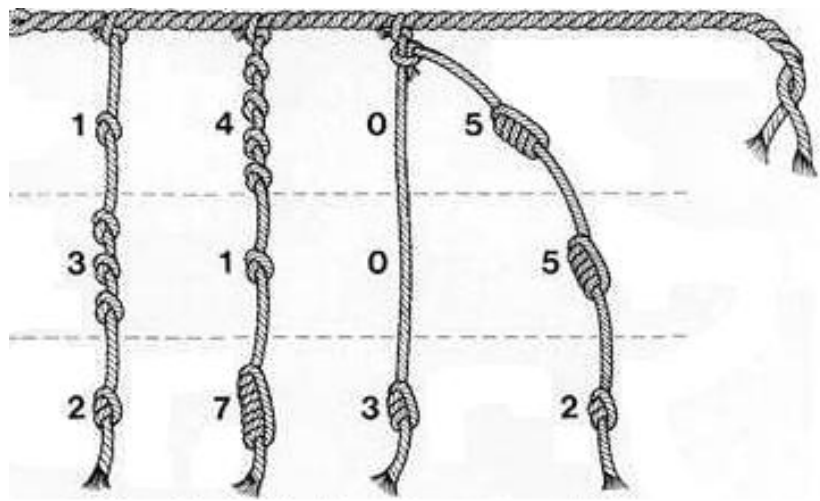

Fig.4 Example of reading numbers in Khipu [8]

In the above example we can see that in the first pendant as we start reading from the end position we encounter '2' which is used to encode using a long knot of 2 turns and hence it is the units place. Next we have 3 and 1 and the number on that pendant becomes 'hundred and thirty-two'. In this way all the values of all chords are added up to get the number encoded on that khipu.

## Characteristics of the Khipu

For storing information on a khipu, its characteristics play a very important role. We will look into the basic characteristics of a khipu which are as follows [8]:

1. Material used for construction
2. Method of twisting of threads and fibers
3. Method of attachment of the pendant to the primary chord
4. Knot directionality
5. Color of the pendants

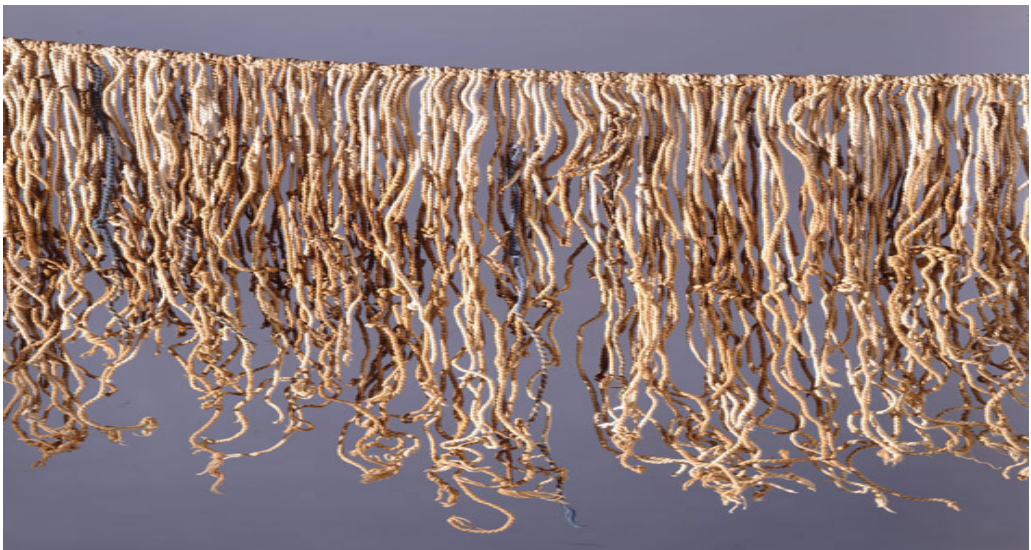

Fig.5 Another example of a khipu [1]

## Material of Construction of the Khipu

Most of the extant khipu are made of cotton fiber or wool. Some of the khipu are also made of llama or alpaca fibers while some are also constructed with the mixture of both llama and alpaca. Few of the khipu are made of vegetable fibers as well. Some khipu are found to be made of human hair.

## Method of Twisting of Fibers

The khipu are formed by twisting together several threads which are in turn formed by twisting together several fibers (Figure 6). Twisting of several fibers is called as 'spinning' and twisting of several threads together is called as 'plying'. Spinning and plying is always done opposite to each other i.e. if spinning is done in clockwise direction then plying is done in anti-clockwise direction [1].

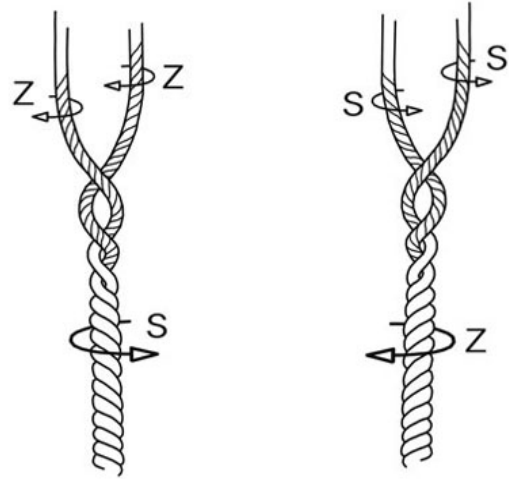

Fig.6 Spinning and plying [8]

**Method of attachment of pendants to the main chord**
Pendants are attached to the primary chord by forming a hitch knot around the main chord. The knot is formed is similar to the tags that we attach to our luggage bags. The hitch knot forms a loop through which the main chord passes. The method of forming the knot is of two types:
  Recto- This type of method forms a loop which is in the front of the main chord.
  Verso –This method forms a loop which is at the back of the main chord.

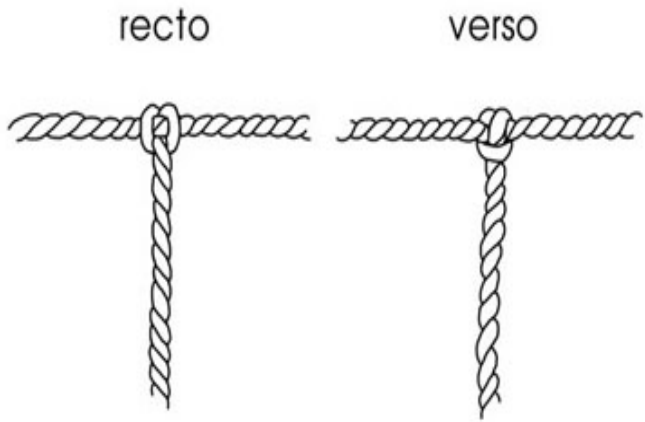

Fig.7 Method of attachment of pendants to main chord [8]

**Knot Directionality**
There are three types of knots Long knot, simple knot and Figure eight knot. All these knots are formed by twisting the threads in one of the two directions. While twisting there forms an axis which slants in either direction. This direction determines the knot directionality. As the directions for spinning and twisting are different, similarly knot orientations are called S and Z [1].

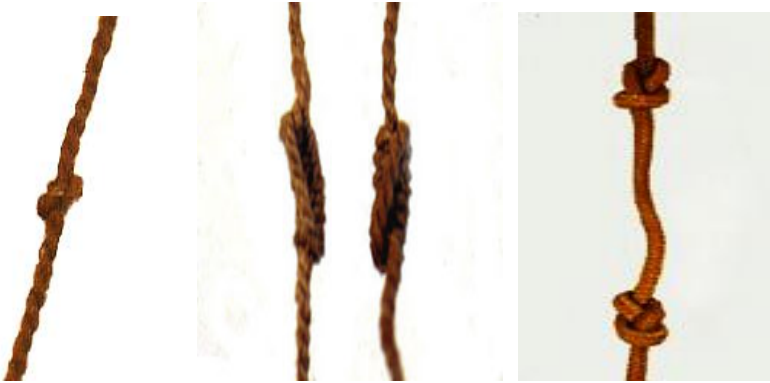

Fig.8 Knot Directionality [1]

**Color of the Pendants of the Khipu**

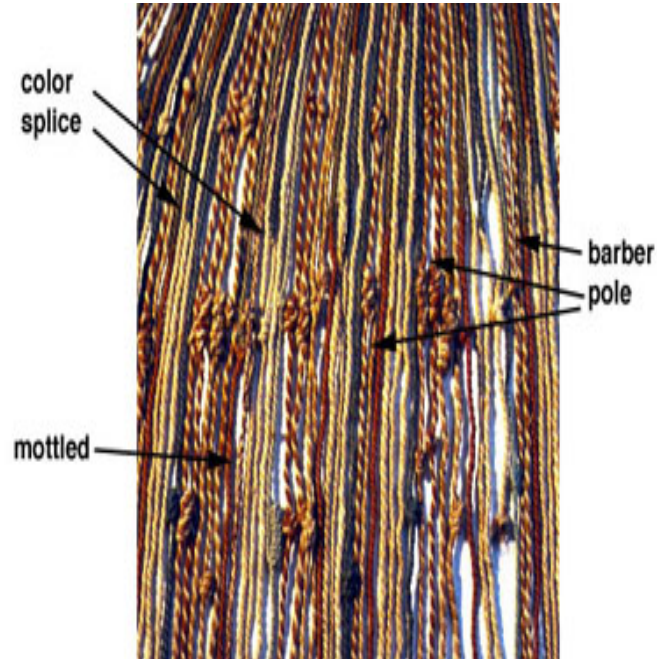

Fig.9 Coloring in Khipu [8]

There are various colors of pendants such as brown, green, blue etc. Predominantly brown color pendants are found but all colors are not solid. Depending on color variation there are two types:

1. Barber Pole- Two threads of different color are plied together to get Barber Pole.
2. Mottled Chord- Spinning threads with two different colors form mottled chord.

**Binary Coding in the Khipu**

Binary is the language of computer which comprehends only the language of 0 and 1. Hence binary language has two choices as shown in Figure 10.

Considering the characteristics of the khipu, most of them have binary choices [27]:

1. Considering the material used: cotton or wool.
2. Color: brown shade and other colors.
3. Spinning and plying can take place in either of two ways, spinning/plying or plying/spinning.
4. The attachment of the pendant chords: recto and verso.
5. Knot directionality is of two types: S or Z.
7. Number class can be decimal or non-decimal.

Considering the above points Gary Urton proposes that seven types of information are coded in binary bits. The knot in the figure above on the right-hand side shows how the 7-bit code can be used. Urton believes that "the binary coding of khipu constituted a means for encoding paired elements that were in relationships of binary opposition to each other, and that, at a semantic level, these relations were of a character known in the literature as markedness relations." [27]

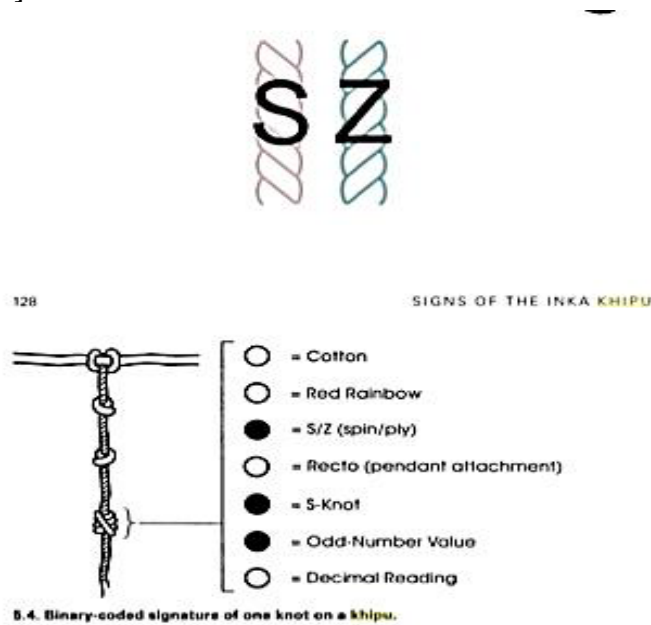

Figure 10. Binary encoding on the khipu [27]

## Conclusions

The paper summarizes some main physical characteristics of the khipu. There is concerted research going on for decoding the khipu enumeration system. There is expectation that the numeration system was devised so that it could also be used for coding narrative. How that might have been done is not clear at this time. The problem of khipu coding is of interest not only for anthropological reasons, but also for abstract mathematical reasons related to information coding.